\documentclass[a4paper]{article}
\title{Maximal cocliques of a strongly regular graph with parameters
 (2048,276,44,36)}
\author{Thomas Jenrich}
\date{2023-03-13}
\begin{document}
\maketitle

\section{Abstract}

This article considers a strongly regular graph with parameters
(2048,276,44,36) that is related to the extended binary Golay code.
That graph is described in \cite{Hub}, \cite{BvL} and \cite{BvM}.

The source package of this article contains a data file that encodes a
sequence of maximal independence sets of that graph, covering
all sizes from 20 to 67 and the size 72, and a Pascal program to check
this assertion and to optionally generate a text file (to be read by the
computer algebra system GAP) that contains the adjacency lists of that
graph and the list of the independence sets.

\section{Introduction}

\subsection{Strongly regular graphs}

A simple loopless finite undirected graph is called strongly regular with
parameter set $(v,k,\lambda,\mu)$, or shortly a srg$(v,k,\lambda,\mu)$,
iff it has exactly $v$ vertices, each of them has exactly $k$ neighbours,
and the number of common neighbours of any two different vertices is
$\lambda$ if they are neighbours and $\mu$ otherwise.

\subsection{The considered srg$(2048,276,44,36)$}

Seemingly, the first published description of the graph is in \cite{Hub},
identified by S.24, in the ``Notes added in proof'', there attributed to
J. H. Conway and M. S. Smith.
But herein, the following construction given in \cite{BvL} will be used:
 ``Take the $2^{11}=2048$ cosets of even weight of the extended binary
 Golay code as vertices, and join two cosets when they have
 representations differing by a vector of weight two.''

Let $X$ the 24-dimensional vector space over the finite field $GF(2)$,
$M$ the set of even-weight elements of $X$, $C$ the extended binary Golay
code, and $R$ the set of those elements of $M$ that have the weight 0,
have the weight 2, or have the weight 4 and the value 1 at the last
position.
Let $w$ denote the weight function, $m$ denote the function that returns
the result of component-wise multiplication of two vectors from $X$,
and $\mathbf{0}$ denote the zero vector of $X$.

If $x \in X$ then $x+x=\mathbf{0}$. If $x,y \in X$ then $x+y=x-y$.

$C$ is known to be a 12-dimensional linear subspace of $M$, consisting of
1, 759, 2576, 759 and 1, resp., vectors of weight 0, 8, 12, 16 and 24,
resp.

$|X|=2^{24}$, $|M|=2^{23}$, $|C|=2^{12}=4096$,
 $|R|={24\choose 0}+{24\choose 2} + {23 \choose 3} = 1+276+1771=2048=2^{11}$.

Let $ x, y \in M$. $w(x-y)=w(x+y)$ equals $w(x)+w(y)-2 \times w(m(x,y))$
and is therefore always even. Thus, $x+y \in M$.

$x$ and $y$ are in the same coset (with respect to $C$) if and only if
$x-y=x+y \in C$.

Let $ x, y \in R \land x \ne y \land z=x+y=x-y$. Because $x \ne y$,
$z \ne \mathbf{0}$. On the other hand, the weight of $z$ is at most 6.
This is obvious if the weight of $x$ or of $y$ is 2 or 0. In
the remaining case, both vectors have weight 4 but also have the value 1
at the last position and so $z$ can have the value 1 in at most $3+3=6$
coordinates.

Thus, $z$ is not in $C$, $x$ and $y$ are not in the same coset. Because
of the cardinalities of $M$, $C$ and $R$, each coset is represented by a
unique element of $R$.

To decide whether the two cosets represented by $x$ and $y$, resp., are
to join, as the above construction says, we have to check whether there
is a $c \in C$ such that the $w((x-y)+c)=2$.

Again, let $z=x+y=x-y$. As just explained, $w(z) \leq 6$. If $w(c)>10$
then $w(z+c) > 10-6=4$. Thus, $w(c)=8 \lor c=\mathbf{0}$ is a derived
precondition.

Case $w(z)=0$:

$z=\mathbf{0}$, $z+c=c$, $w(z+c)=w(c) \ne 2$. The two cosets are not to
join.

Case $w(z)=2$:

If $c=\mathbf{0}$ then $w(z+c)=w(z)=2$. The two cosets are to join.

Case $w(z)=4$:

If $c=\mathbf{0}$ then $w(z+c)=w(z)=4>2$.
If $w(c)=8$ then $w(z+c) \geq 8-4=4>2$.
The two cosets are not to join.

Case $w(z)=6$:

If $c=\mathbf{0}$ then $w(z+c)=w(z)=6>2$. It remains to scan
 $\{ c \in C : w(c)=8 \}$.

Because $w(c)=w(z)+2$, $w(z+c) \geq 2$.

Clearly, if $w(z+c)=2$ then the two cosets are to join.

If $w(z+c)=2$ then all 6 non-zero positions of $z$ are non-zero positions
of $c$.

Assume that there is a $d \in C$ such that $w(d)=8$ and $w(z+d)=4$.
Then 5 of the non-zero positions of $d$ are non-zero positions of $z$
and thus non-zero positions of $c$. This implies $w(c+d) \leq 6$ and
finally $c=d$. But in that case $w(z+d)=w(z+c)=2$.

So, if $w(z+c)=4$ then the two cosets are not to join.

If $w(z+c)>4$, we cannot immediately decide whether the cosets are to join,
continue the search if not all $c \in C$ with weight 8 were checked.

\subsection{Cliques and cocliques}

The parameters of the considered graph imply that the independence
number (size of the largest coclique) cannot exceed 85: The smallest
eigenvalue of the $\{0,1\}$-adjacency matrix of the graph is $-12$. Thus,
the Delsarte bound for the maximum coclique size is
$2048/(1+276/12) = 2048/24 = 85 + 1/3 $.

The current survey \cite{BvM} contains in particular many propositions on
maximal cliques and cocliques of individual strongly regular graphs.
In the case of the graph considered herein, the subsection on cliques is
rather precise but the subsection on cocliques just stated that the
independence number is in the range from 50 to 84, giving neither a proof
nor a dedicated reference.

By an (incomplete) extensive search I have found maximal (i.e., not
extensible) cocliques of any size from 20 to 67 and of size 72.

The source package of this article includes the files SRG2048C2.PAS and
SRG2048C.DAT .

The Pascal source file SRG2048C2.PAS has been (and can be) used to check
that the vertex subsets encoded in SRG2048C.DAT represent one maximal
coclique of each size from 20 to 67 and two maximal cocliques of size 72,
and that the subgraphs induced by the complements of the two vertex subsets
of size 72 are not isomorphic (to each other).

\section{Computations}

SRG2048C2.PAS is a console program for compilers compatible with Borland's
Turbo Pascal 4.0.

The binary vectors of length 24 are encoded as (binary) integer values in
the range from 0 to $2^{24}-1$ (type name \texttt{t\_bin24}). The addition
and subtraction of those vectors is performed by applying the binary
operator \texttt{xor} to the encoding integers.

The generator matrix for the Golay code given by the
constant \texttt{golay\_bases\_str} as an array of 12 character strings
of length 24 has been taken from the respective chapter of \cite{BvM}.

In the initial stage, the variable \texttt{golay8\_code} is filled
with the 759 code words of weight 8 (routine \texttt{init\_golay8\_code})
and the variable \texttt{rep} is filled with the representations of the
2048 cosets (routine \texttt{init\_rep}).

The function routine \texttt{adja} calculates whether the two cosets
whose representations are given as parameters are adjacent. For this
purpose, it calls the function routine \texttt{minw}. The task of
\texttt{minw} is to find a weight 8 code word such that the
weight of the sum of that code word and the vector encoded by the
function parameter value is at most 4 or as small as possible.
If \texttt{minw} returns a value that is not 0, 2 or 4, then
\texttt{adja} shows that value together with the message
``invalid distance'' and stops the execution. Probably, that can not
happen, but I do not have a proof.

The routine \texttt{read\_and\_check\_cocliques} expects that the binary
file SRG2048C.DAT is in the current working directory and that it contains
a sequence of vector sets in the form of strings where each string starts
with a byte giving the size of the then following sequence of unsigned
three-byte integer values (that encode the vectors).
If the size is not in the range from 2 to 85, or if one of the vectors
is not a proper coset representation, a message is displayed and the
execution stops.

For each vector sequence found in the file, its size is displayed and it
is checked to establish a maximal coclique of that size.

The size of the largest cocliques that I have found is 72. In order to
show that there are (at least) two structurally different cocliques of
that probably maximal size, the file contains two vertex sets of size 72
and the program calls the routine \texttt{check\_external\_relation} for
each coclique of size at least 72 found in the file.

The first part of that routine counts the vertices out of the considered
set having certain numbers of neighbours in that set and
shows the resulting statistic. Because the result was the same
(namely \texttt{8:480 10:960 12:536} ) for all cocliques of size 72 that
I've investigated, a more detailed analysis has been added: The number of two-sets of vertices in the considered set having no common neighbour that is adjacent to exactly 8 vertices in the considered set is calculated and displayed.
The results for the two vertex sets of size 72 contained in SRG2048C.DAT
were 166 and 276, resp.

\subsection{Optional data export for GAP}

SRG2048C2.PAS has been build from its predecessor SRG2048C.PAS (published
with version 1 of this article) by adding optional instructions to
generate a text file SRG2048C.g (to be read by the computer algebra
system GAP; about 2.6 MB). These instructions will be compiled into the
executable file only if the compiler symbol GAP is defined.
The content of SRG2048C.g would look like

\ttfamily

\ A:= ..... ;

\ MIS:= ..... ;

\ LoadPackage("grape");;

\ Gra:=Graph(Group(()),[1..2048],OnPoints,

\ function(x,y) return (x in A[y]); end, true);

\rmfamily
\smallskip

where the two five-dot sequences respectively stand for the list of
adjacency lists of the vertices of the considered graph and the list
of the maximal independence sets represented by lists of vertex
numbers.

\subsection{Compilation and execution}

The command line instances (TPC, DCC32 and fpc, resp.) of the compilers
have been used. In the case of fpc, the command line option -Mtp has been
given (for compatibility with Turbo Pascal).

\subsubsection{Without the optional data export}

Actual compilations and executions have been done on two different
computer systems:

System 1:

1 GHz Intel Pentium(R) III, running MS Windows 98SE.
Compilers and execution times:

Turbo Pascal 5.5 : 2:12.81 min

Turbo Pascal 7.01 : 1:43.12 min

Delphi 4.0 build 5.37 : 5.00 s

Free Pascal 2.4.4 for i386 : 6.75 s

System 2:

2.8 GHz Intel Pentium(R) Dual-Core E5500, running (Linux distribution)
Lubuntu 20.04 (64 bit).

Compiler and execution time: Free Pascal 3.0.4 for x86\_64 : 2.43 s

\subsubsection{With the optional data export (compiler symbol GAP defined)}

The command line option -dGAP should work for all mentioned compilers.

\smallskip

Actual compilations and executions have been done on system 2 as
described above.

Compiler and execution time: Free Pascal 3.0.4 for x86\_64 : 10.52 s

\section{Remark}

In August 2021, a few weeks after the appearance of version 1 of this
article, Ivan Mogilnykh wrote to me that some time before he had
together with his colleague Denis Krotov found two maximal
independence sets of size 72 that are structurally different (from
each other). He talked about (but did not publish) the first
of those two sets at the Mal'tsev meeting 2020 (abstract in \cite{Mog}).
According to the eMails from I. Mogilnykh and D. Krotov that I received
today, the calculation of the invariant used in my program to
structurally distinguish the two sets of size 72 (giving 166 and 276,
resp.) done for the sets found by them resulted in 336 and 166, resp.
Thus, at least the first of their sets is structurally different from
the two sets that I found.

\vspace{0.1in}

Author's eMail address: thomas.jenrich@gmx.de

\end{document}